\documentclass[a4paper,12pt]{amsart}
\usepackage{amsmath,amssymb}
\usepackage[shortlabels]{enumitem}
\usepackage[cp1251]{inputenc}

\makeatletter
\@namedef{subjclassname@2020}{%
\textup{2020} Mathematics Subject Classification}
\makeatother

\newtheorem{thm}{Theorem}

\newenvironment{proof*}{\vskip 2mm\noindent {}}{\hfill $\Box$ \vskip 2mm}

\title[A Picard little theorem for entire functions of matrices]
{A Picard little theorem for entire functions of matrices}

\author{Oleg Mushkarov}
\address{Institute of Mathematics and Informatics, Bulgarian Academy of Sciences, Acad. G. Bonchev 8,
1113 Sofia, Bulgaria}
\email{muskarov@math.bas.bg}

\author{Nikolai Nikolov}
\address{Institute of Mathematics and Informatics, Bulgarian Academy of Sciences, Acad. G. Bonchev 8,
1113 Sofia, Bulgaria
\vspace{1mm}
\newline Faculty of Information Sciences, State University of Library Studies and Information
Technologies, Shipchenski prohod 69A, 1574 Sofia, Bulgaria}
\email{nik@math.bas.bg}

\thanks{The second named author was partially supported by the Bulgarian National Science Fund,
Ministry of Education and Science of Bulgaria under contract KP-06-N82/6.}


\begin{document}


\maketitle

\vspace{-8mm}
\begin{abstract} An analog of Picard's little theorem for entire functions of
matrices is proved.
\end{abstract}

\section{Introduction}

The Picard little theorem says that any non-constant entire
function omits at most one complex value, i.e. its range is either
the whole complex plane or the complex plane minus a single point.
This theorem together with Picard's great theorem have a long
history and many important generalizations that have been the
impetus for powerful developments in complex analysis such as the
Nevanlinna theory for value distribution of meromorphic functions
(see e.g. \cite{Ha}).

The purpose of this note is to prove an analog of Picard's little
theorem for entire functions of complex matrices. More precisely,
let $M_n(\mathbb{C})$ be the ring of all $n\times n$ matrices with
complex entries. Then any entire function $f:\mathbb{C}\to
\mathbb{C}$ with power series $$f(z)=\sum_{k=0}^\infty a_kz^k$$
defines a map $f : M_n(\mathbb{C}) \to M_n(\mathbb{C})$ by
setting
$$ f(A)=\sum_{k=0}^\infty a_kA^k$$
for any $A\in M_n(\mathbb{C})$. By analogy with Picard's little
theorem one may ask what can be said about the range of the map $f
: M_n(\mathbb{C}) \to M_n(\mathbb{C})$. This question has been
already considered in \cite [Problem 210*]{PS}, and in \cite{Mo} a
criterion for surjectivity of $f$ has been proved for matrices
with entries in an arbitrary algebraically closed field .

The main purpose of this note is to show that the range of $f$ can
be described almost completely by means of the {\it totally
ramified values} of $f$, i.e. the complex numbers $a$ such that
all roots of the equation $f(z)=a$ have multiplicity at least 2.

For $a\in \mathbb{C},$ let $E_a$ be the set of all matrices $A\in
M_n(\mathbb{C})$ which have an eigenvalue $a,$ and let $S_a\subset
E_a$ be the set of matrices whose Jordan forms have at least one
non-trivial Jordan block corresponding to $a$ (i.e. the eigenspace
of $A$ corresponding to $a$ has non-zero codimension). Our analog
of Picard's little theorem for entire functions of matrices is the
following

\begin{thm}\label{rv} Let $f:\mathbb{C} \to\mathbb{C}$ be a non-constant
entire function.

(i) If $f$ omits a complex value $a$, then
$f(M_n(\mathbb{C}))= M_n(\mathbb{C})\setminus E_a$.

(ii) If $f$ takes on all complex values and has no totally
ramified values, then $f(M_n(\mathbb{C}))= M_n(\mathbb{C})$.

(iii) If $f$ has one totally ramified value $a$, then
$f(M_n(\mathbb{C}))= M_n(\mathbb{C})\setminus S^f_a$, where
$\emptyset \neq S^f_a\subset S_a$. Moreover $S^f_a = S_a$ if and
only if each root of the equation $f(z)=a$ has a multiplicity at
least $n$. In particular, $S^f_a=S_a$ for $n=2$.

(iv) If $f$ has two totally ramified values $a$ and $b$, then
$f(M_n(\mathbb{C}))= M_n(\mathbb{C})\setminus (S^f_a\bigcup
S^f_b)$. Moreover $S^f_a$ and $S^f_b$ are non-empty proper subsets
of $S_a$ and $S_b$ for $n\geq 3$, while $S^f_a=S_a$ and $S^f_b=S_b$ for
$n=2$.
\end{thm}

To prove Theorem \ref{rv} we use some well-known algebraic facts
about square complex matrices, Picard's little and great theorems
as well as  the well-known fact  that any non-constant entire
function has at most two totally ramified values.

The needed auxiliary facts for square complex matrices and totally
ramified values of entire functions are given in Section 2 and
then in Section 3 we prove Theorem 1. Finally, in Section 4 we
discuss some examples of entire functions satisfying the
conditions of the four options listed in Theorem 1.
\smallskip

\textbf{Acknowledgment.} We wish to thank the anonymous referee
for pointing to us that parts (i) and (ii) of Theorem \ref{rv} are
already proven in \cite{Mo}, a reference we were not aware of
while writing the present note.

\section{Auxiliary facts}

We first recall some well-known algebraic facts about square
complex matrices.

A Jordan block of size $k$ is a $k\times k$ matrix of the form

\begin{displaymath}
J_k(\lambda)=
\begin{pmatrix}
\lambda & 1     &  \\
  & \lambda & 1 &  \\
  & &\ddots
            &1 \\
  &       &        &\lambda\\
\end{pmatrix}
\end{displaymath}
where the missing entries are all zero. Every square complex
matrix $A$ is similar to a block diagonal matrix

\begin{displaymath}
J_A=
\begin{pmatrix}
A_1      &  \\
  &A_2  &  \\
  & &\ddots
             \\
  &       &        &A_m\\
\end{pmatrix}
\end{displaymath}
where each block $A_s$ is a Jordan block. The matrix $J_A$ is
unique up to the order of the Jordan blocks and is called the
Jordan form of $A$. The entries on its main diagonal are equal to
the eigenvalues of $A$.

Let $f:\mathbb{C} \to\mathbb{C}$ be an entire function.
\smallskip

\textbf{Fact 1.}\label{f1} For any  $X\in M_n(\mathbb{C})$ and
any non-singular $P\in M_n(\mathbb{C})$ we have $f(P^{-1}XP)=P^{-1}f(X)P .$
\smallskip

\textbf{Fact 2.}\label{f2} For any Jordan block $J_k(a)$ the
following identity holds true

\begin{displaymath}
f(J_k(a))=
\begin{pmatrix}
f(a) & f'(a) & \dots & \dots & \frac{f^{(k-1)}(a)}{(k-1)!} \\
& f(a) & f'(a) & \dots & \frac{f^{(k-2)}(a)}{(k-2)!} \\
        &           &\ddots &\ddots &\vdots\\

  &            &            &f(a) &f'(a)\\
    &          &      & &f(a) \\
\end{pmatrix}
\end{displaymath}
where the missing entries are all zero.

\smallskip

\textbf{Fact 3.}\label{f3} Let $A$ be a block diagonal matrix with
diagonal blocks $A_1,A_2,\dots,A_m$. Then $f(A)$ is a block
diagonal matrix with diagonal blocks $f(A_1),f(A_2),\dots,
f(A_m).$
\smallskip

{\it Proof of Facts 1,2, and 3.} It is enough to prove each of these
facts for the function $f(z)=z^m$, where $m$ is a positive
integer. This can be done easily by induction on $m$.
\smallskip

Given an upper triangular matrix $A\in M_n(\mathbb{C})$ denote by
$d_0(A)$ its main diagonal and by $d_k(A),$ $1\leq k\leq n-1,$ the
$k$-th diagonal of $A$ above $d_0(A)$.
\smallskip

\textbf{Fact 4.}\label{f4} Let $A\in M_n(\mathbb{C})$ be an upper
diagonal matrix such that all entries on $d_0(A)$ are $0$ and all
entries on $d_1(A)$ are equal to $a$. Then for any $1\leq k\leq
n-1$ all entries on $d_1(A^k),\dots ,d_{k-1}(A^k)$ are $0$ and all
entries on $d_k(A^{k})$ are equal to $a^k$.

{\it Proof.} Induction on $k$.
\smallskip

\textbf{Fact 5.}\label{f5} Let $A\in M_n(\mathbb{C})$ be an upper
triangular matrix whose entries on $d_1(A)$ are all equal to
$a$.Then the Jordan form of $A$ has a single Jordan block
$J_n(\lambda)$ if and only if $a\neq 0$ and all entries on $d_0(A)$
are equal to $\lambda$.
\smallskip

{\it Proof.} Note first that all entries on $d_0(A)$ are equal to
$\lambda$ since the similar matrices have equal eigenvalues. Hence
the Jordan form of $A$ consists only of Jordan blocks of the form
$J_s(\lambda)$. Denote by $k$ the largest size of such a Jordan
block. Then $k=n$ if and only if $(A-\lambda I)^{n-1}\neq 0$. On
the other hand, we know from Fact 4 that the entry on
$d_{n-1}((A-\lambda I)^{n-1})$ is equal to $a^{n-1}$ which proves
Fact 5.  \qed
\smallskip

We will need also the following analytic facts.
\smallskip

\textbf{Fact 6.} \cite[p. 45]{He} Any entire function $f$ has at
most two totally ramified values. If they are two, say $a$ and
$b$, then each of the equations $f(z)=a$ and $f(z)=b$ has a root
of multiplicity $2$.
\smallskip

\textbf{Fact 7.} If a non-constant entire function omits a value,
then it has no totally ramified values.
\smallskip

{\it Proof.} See Example 3 for a more general result.

\section{Proof of Theorem \ref{rv}}

Given a matrix $A\in M_n(\mathbb{C})$ denote by $J_A$ its Jordan
form. Then it follows from Fact 1 that the equation $f(X)=A$ has a
solution if and only if there is a matrix $X\in M_n(\mathbb{C}) $
such that the matrices $f(J_X)$ and $J_A$ are similar. This
together with Facts 3, 2, 5, and 7 implies that the
equation $f(X)=A$ has a solution if and only if for any eigenvalue
$a$ of $A$ there exists $z_a$ such that $f(z_a)=a$ and $f'(z_a)\neq 0$.
This proves statements (i) and (ii).

Suppose now that $f$ has one totally ramified value $a$. If $a$ is
not an eigenvalue of $A$ then the equation $f(X)=A$ has a solution
since otherwise we will obtain another totally ramified value of
$f(z)$, a contradiction. The reasoning above shows also that if
$a$ is an eigenvalue of $A$ and $J_n(a)\in S_a$ then the equation
$f(X)=J_n(a)$ has no solutions, i.e.
$f(M_n(\mathbb{C}))=M_n(\mathbb{C})\setminus S^f_a,$ where
$\emptyset \neq S^f_a\subseteq S_a.$ Moreover, if the equation
$f(z)=a$ has a root $z_a$ with multiplicity less than $n$, then
$S^f_a$ is a proper subset of $S_a$. Indeed, it follows from Fact
2 that the Jordan form of $f(J_{n}(z_a))$ has at least one
non-trivial Jordan block which shows that $f(J_{n}(z_a)\in S_a$.
To prove (iii), it remains to show that if each root of the
equation $f(z)=a$ has a multiplicity at least $n$, then
$S^f_a=S_a$. Suppose that $f(X)=A$, where $A\in S_a$. Then the
Jordan form of $A$ has at least one non-trivial Jordan block
$J_k(a)$. On the other hand, Fact 2 shows that each of these
Jordan blocks is trivial which is a contradiction.

Finally, suppose that $f$ has two totally ramified values $a$ and
$b$. Then by Fact 6 any of the equations $f(z)=a$ and $f(z)=b$ has
a root of multiplicity 2 and (iv) follows by using the same
reasoning as in the proof of (iii). \qed
\smallskip

\textbf{Remark.} Denote by $GL_n(\mathbb{C})$ the group of
non-singular $n\times n$ matrices with complex entries. Then any
holomorphic function $f:\mathbb{C}\setminus\{0\} \to\mathbb{C}$
with Laurent series
$$f(z)=\sum_{k=-\infty}^\infty a_kz^k$$
defines a map $f : GL_n(\mathbb{C}) \to M_n(\mathbb{C})$ by
setting
$$ f(A)=\sum_{k=-\infty}^\infty a_kA^k .$$
Since $g(z)=f(e^z)$ is an entire function, it follows that $f$
omits at most one value as well as that $f$ has at most two
ramified values. Note also that $f$ has no ramified values if it
omits a value. Hence one can repeat the proof of Theorem \ref{rv}
line by line and see that its statement is true for holomorphic
functions on $\mathbb{C}\setminus\{0\}$ as well. This follows also
as a corollary of Theorem \ref{rv} since
$f(GL_n(\mathbb{C}))=g(M_n(\mathbb{C}))$.

\section{Examples}

In this section we discuss some examples of entire functions
satisfying the conditions of the four options listed in Theorem
\ref{rv}.
\smallskip

\textbf{Example 1.} (i) Any quadratic polynomial has one totally
ramified value.

(ii) The polynomial $P(z)=z^k(z-1)$ ($k\ge 2)$ has
no totally ramified values.

(iii) The general fact for polynomials is that they have at most
one totally ramified value. This follows by the simple observation
that if $a$ is such a value for a polynomial $P$, then the zeros
of $P-a$ contribute at least
$\displaystyle\frac{\mbox{deg}(P)}{2}$ to the zeros of $P'$
counted with multiplicity.
\smallskip

\textbf{Example 2.} Consider the function
$$f(z)=\displaystyle\frac{a-b}{2}\sin(cz+d)+\frac{a+b}{2},$$
where $a\neq b,c\neq 0,d\in \mathbb{C}.$ Then
$f'(z)=\displaystyle\frac{(a-b)c}{2}\cos(cz+d)=0 \Leftrightarrow
\sin(cz+d)=\pm 1$ and we conclude that $f$ has two totally
ramified values $a$ and $b$. Moreover, by \cite[Theorem 2.1]{BQL},
these are all entire functions of order at most 1 (i.e.
$f(z)=\mbox{O}(e^{|z|^\gamma})$ for any $\gamma>1$) with this
property.
\smallskip

\textbf{Example 3.} As we know by Picard's great theorem, for any
transcendental entire function $f$ (i.e. not a polynomial) there
is at most one $a\in\mathbb{C}$ such that the equation $f(z)=a$
has finitely many roots. In this case $f=a+Pe^g$, where $P$ is a
monic polynomial and $g$ is a non-constant entire function. Then
it follows by Nevanlinna's theorem on deficient values
\cite[Theorem 2.4]{Ha} that $a$ is the only possible totally
ramified value of $f$ and this is so if and only if all zeros of
$P$ have multiplicity greater than $1$. In particular, if a
non-constant entire function $f$ omits a value $a$, then $P=1$ and
$f$ has no totally ramified values (Fact 7).

\end{document}